\documentclass[11pt,a4paper]{elsarticle}

\usepackage{amssymb}
\usepackage[ansinew]{inputenc}
\usepackage{amsmath}
\usepackage{amsthm}

\usepackage{footmisc}
\usepackage{mathrsfs} %ÓÃÀ´´ò»¨ÌåµÄÓ¢ÎÄ×ÖÄ¸
\usepackage{latexsym, enumerate}
\usepackage{amsfonts}

\newtheorem{theorem}{\textbf{Theorem}}[section]
\newtheorem{lemma}{\textbf{Lemma}}[section]
\newtheorem{proposition}{\textbf{Proposition}}[section]
\newtheorem{corollary}{\textbf{Corollary}}[section]
\newtheorem{remark}{\textbf{Remark}}[section]
\newtheorem{definition}{\textbf{Definition}}[section]

\makeatletter  %´òÂÞÂíÐòºÅ
    
    \newcommand{\Rmnum}[1]{\expandafter\@slowromancap\romannumeral #1@}
  \makeatother

\newcommand{\bm}[1]{\mbox{\boldmath{$#1$}}}  %ÂÞÂí×Ö·ûºÚÌå

\allowdisplaybreaks[4]

\def\be{\begin{equation}}
\def\ee{\end{equation}}
\def\ube{\begin{equation*}}
\def\uee{\end{equation*}}
\def\ube{\begin{equation*}}
\def\uee{\end{equation*}}
\def\bee{\begin{eqnarray}}
\def\beee{\begin{eqnarray*}}
\def\eee{\end{eqnarray}}
\def\eeee{\end{eqnarray*}}
\def\bt{\begin{theorem}}
\def\et{\end{theorem}}
\def\bl{\begin{lemma}}
\def\el{\end{lemma}}
\def\br{\begin{remark}}
\def\er{\end{remark}}
\def\bp{\begin{proposition}}
\def\ep{\end{proposition}}
\def\bc{\begin{corollary}}
\def\ec{\end{corollary}}
\def\bd{\begin{definition}}
\def\ed{\end{definition}}
\def\non{\nonumber}

\journal{Journal of Mathematical Analysis and Applications}

\begin{document}

\begin{frontmatter}



\title{Vanishing viscosity limit for a coupled Navier-Stokes/Allen-Cahn system}

 \author[label1]{ Liyun Zhao\footnote{Corresponding author: Haiyang Huang.\\
  E-mail addresses: zhaoliyun@mail.bnu.edu.cn(L. Zhao), gbl@mail.iapcm.ac.cn(B. Guo), hhywsg@bnu.edu.cn(H. Huang).}}
  \address[label1]{ School
of Mathematical Sciences, and Key Laboratory of Mathematics and
Complex Systems (Ministry of Education), Beijing Normal University,
Beijing 100875, China.}
 \author[label2]{ Boling Guo}
 \address[label2]{Nonlinear Center
for Studies, Institute of Applied Physics and Computational
Mathematics, PO Box 8009, Beijing 100088, China. }
 \author[label1]{ Haiyang Huang}

\begin{abstract}
In this paper, we study the vanishing viscosity limit for a coupled
Navier-Stokes/Allen-Cahn system in a bounded domain. We first show
the local existence of smooth solutions of the Euler/Allen-Cahn
equations by modified Galerkin method. Then using the boundary layer
function to deal with the mismatch of the boundary conditions
between Navier-Stokes and Euler equations, and assuming that the
energy dissipation for Navier-Stokes equation in the boundary layer
goes to zero as the viscosity tends to zero, we prove that the
solutions of the Navier-Stokes/Allen-Cahn system converge to that of
the Euler/Allen-Cahn system in a proper small time interval. In
addition, for strong solutions of the Navier-Stokes/Allen-Cahn
system in 2D, the convergence rate is $c\nu^{1/2}$.
\end{abstract}


\end{frontmatter}
%%%%%%%%%%%%%%%%%%%

\section{Introduction}
In this paper, we are concerned with the vanishing viscosity limit
for the following Navier-Stokes/Allen-Cahn system in
$\Omega\times(0,+\infty)$: \bee && {\bf u}_t+({\bf
u}\cdot\nabla){\bf u}+\nabla P_1 = \nu\Delta{\bf u}-\lambda\nabla\cdot(\nabla\phi\otimes\nabla\phi),\label{Navier-Stokes} \\
&&\nabla\cdot {\bf u}= 0,\label{divergence free}\\
&&\phi_t+({\bf u}\cdot\nabla)\phi =
\gamma(\Delta\phi-f(\phi)),\label{Allen Cahn} \eee with initial data
 \be
 {\bf u}(x,0)={\bf u}_0(x),\quad \phi(x,0)=\phi_0(x),\ \ x\in\Omega,\label{4}
 \ee
and boundary conditions:
 \be
 {\bf u}(x,t)=0,\quad \partial_\textbf{n}\phi(x,t)=0, \ \ (x,t)\in \partial\Omega\times(0,+\infty).\label{BC1}
 \ee
 Here, $\Omega\subset{\mathbb{R}^n}$
($n=2,3$) is a bounded domain with smooth boundary $\partial\Omega$.

System \eqref{Navier-Stokes}-\eqref{BC1} can be viewed as a phase
field model, which describes the motion of a mixture of two
incompressible viscous fluids with the same density and viscosity
(see \cite{LS03,Yue04}). The two fluids are macroscopically
immiscible and separated by a thin interface. \eqref{Navier-Stokes}
is the linear momentum equation, where ${\bf u}$ is the velocity
field of the mixture, $\phi$ and $P_1$ denote the phase function and
pressure, respectively. $\nabla\phi\otimes\nabla \phi$ denotes the
induced elastic stress, which is a $n\times n$ matrix whose
$(i,j)$-th entry is $\frac{\partial\phi}{\partial
x_i}\frac{\partial\phi}{\partial x_j}$ for $1\leq i,j\leq n$.
\eqref{divergence free} implies the incompressibility of both fluids
in the mixture and \eqref{Allen Cahn} is the phase equation of
Allen--Cahn type.  $f(\phi) = \frac{1}{\varepsilon^{2}}(\phi^3 -
\phi)$. $\nu$, $\lambda$, $\gamma$, $\varepsilon$ are positive
constants, representing the kinematic viscosity, the surface
tension, the mobility and the width of the interface, respectively.
${\bf u}|_{\partial\Omega}=0$ is non-slip boundary condition. And
$\partial_\textbf{n}\phi|_{\partial\Omega}=0$  means that the
diffused interface perpendicularly contacts the boundary of the
domain, where $\textbf{n}$ is the outward unit normal to the
boundary $\partial\Omega$.

From another point of view, system
\eqref{Navier-Stokes}-\eqref{Allen Cahn} is closely related to
liquid crystal model, Magnetohydrodynamics (MHD) equations, and
viscoelastic system with infinite Weissenberg number, see
\cite{XZL10}.

System \eqref{Navier-Stokes}-\eqref{BC1} has a basic energy law \be
\frac{1}{2}\frac{d}{dt}\left(|{\bf
u}|^2+\lambda|\nabla\phi|^2+\frac{\lambda}{2\varepsilon^2}
|(\phi^2-1)|^2\right)= -(\nu|\nabla{\bf
u}|^2+\lambda\gamma|\Delta\phi-f(\phi)|^2),
 \label{basic energy law}
\ee where $|\cdot|$ denotes the $L^2$-norm. \eqref{basic energy law}
can be derived by multiplying \eqref{Navier-Stokes} by ${\bf u}$,
\eqref{Allen Cahn} by $\lambda(-\Delta\phi+f(\phi))$, then adding
them up and integrating over $\Omega$, also using \eqref{divergence
free} and \eqref{BC1}. Thanks to this important energy law
\eqref{basic energy law}, it can be proved that global weak
solutions of system \eqref{Navier-Stokes}-\eqref{BC1} exist in both
2D and 3D case. With regular initial data, strong solutions exist
globally in 2D, but locally in 3D in a short time interval
$[0,T_\nu)$, where $T_\nu$ depends on $\nu$. Hence we just consider
weak solutions of \eqref{Navier-Stokes}-\eqref{BC1} in the 3D case
in the proof of vanishing viscosity limit.

At the limit case, namely $\nu=0$, Navier-Stokes/Allen-Cahn system
formally becomes the following Euler/Allen-Cahn system:
\bee&&{\bf v}_t+({\bf v}\cdot\nabla){\bf v}+\nabla P_2=-\lambda\nabla\cdot(\nabla\psi\otimes\nabla\psi),  \label{Euler} \\
&&\nabla\cdot {\bf v} = 0,  \\
&&\psi_t+({\bf v}\cdot\nabla)\psi =
\gamma(\Delta\psi-f(\psi)),\label{psi}\eee in
$\Omega\times(0,+\infty)$ with initial data
 \be
 {\bf v}(x,0)={\bf v}_0(x),\quad \psi(x,0)=\psi_0(x),\ \ x\in\Omega,
 \ee
and boundary conditions:
 \be
 {\bf v}\cdot{\bf n}(x,t)=0,\ \  \partial_\textbf{n}\psi(x,t)=0, \ \ (x,t)\in \partial\Omega\times(0,+\infty).\label{BC2}
 \ee
For simplicity, we take ${\bf v}_0(x)={\bf u}_0(x)$ and
$\psi_0(x)=\phi_0(x)$. By taking the limit of weak formulation of
the problem (1.1)-(1.5) as $\nu\rightarrow 0$, we obtain the weak
formulation of (1.7)-(1.11), from which we can derive the system
(1.7)-(1.11). Thus the boundary condition of ${\bf v}$ is ${\bf
v}\cdot{\bf n}|_{\partial\Omega}=0$, which describes that the
boundary is impermeable. The basic energy law of system
\eqref{Euler}-\eqref{BC2} is \be\label{energy-law-2}
\frac{1}{2}\frac{d}{dt}\left(|{\bf
v}|^2+\lambda|\nabla\psi|^2+\frac{\lambda}{2\varepsilon^2}
|(\psi^2-1)|^2\right)= -\lambda\gamma|\Delta\psi-f(\psi)|^2. \ee

The main purpose of this article is to show that the solution of the
viscid system \eqref{Navier-Stokes}-\eqref{BC1} converges to that of
the inviscid system \eqref{Euler}-\eqref{BC2} in a short time period
as the viscosity goes to zero. The inviscid limit helps us to
understand turbulent phenomena governed by the viscous equations,
see Navier-Stokes equations for an example (see e.g.
\cite{PJ1995,Kato1984,Kato-JFA1984,BL-2004,Temam-Wang1997,Temam-Wang2002}).

Let us first recall the classical issue of vanishing viscosity limit
of Navier-Stokes equations in a bounded domain. In order to deal
with the mismatch of the boundary conditions between Navier-Stokes
system and Euler system, Kato in \cite{Kato1984} introduced a
boundary layer function defined on a boundary strip $\Gamma_{c\nu}$
with width $c\nu$, and proved that the solution ${\bf u}^\nu$ of
Navier-Stokes equation converges to that of Euler equation in
$L^2$-norm, provided \be\nu\int^T_0|\nabla{\bf
u}^\nu(t)|_{L^2(\Gamma_{c\nu})}^2dt\rightarrow 0,\ \mbox{as}\
\nu\rightarrow0,\label{Kato}\ee
 where $[0,T]$ is the interval on which the smooth solution of Euler equation
 exists. Later,
Temam and Wang in \cite{Temam-Wang1997} and \cite{Wang01} improved
Kato's condition by replacing the total gradient in \eqref{Kato}
with tangential derivatives only, but required a slightly thicker
boundary layer.

With the same boundary layer function in\cite{Kato1984}, and under
the same assumption \eqref{Kato}, we can also prove that weak
solution of \eqref{Navier-Stokes}-\eqref{BC1} converges to smooth
solution of \eqref{Euler}-\eqref{BC2} in the basic energy space in
$[0,T]$ as the viscosity tends to zero. In particular, for 2D case,
strong solution of \eqref{Navier-Stokes}-\eqref{BC1} converges to
smooth solution of \eqref{Euler}-\eqref{BC2} at the rate
$c\nu^{1/2}$. The main difficulty in the proof is to deal with
coupled terms of ${\bf u}$, ${\bf v}$, $\phi$ and $\psi$. A key
point is that the highest order coupled terms can be canceled, similar to
the derivation of the energy law \eqref{basic energy law}. To
estimate other coupled terms, high regularity of solutions to
\eqref{Euler}-\eqref{BC2} is required. Therefore we first establish
the local existence and uniqueness of smooth solutions to
Euler/Allen-Cahn system \eqref{Euler}-\eqref{BC2}. This proof
requires the energy law \eqref{energy-law-2}, so we employ the
modified Galerkin method introduced in \cite{LL95}, namely, we only
project \eqref{Euler} into finite dimensional space, but solve
\eqref{psi} for any given ${\bf v}$ by the fixed point theorem, such
that the approximate solutions still satisfies the energy law
\eqref{energy-law-2}.

This article is organized as follows. In section 2, we introduce
mathematical preliminaries and state our main results. In section 3,
we establish the local existence and uniqueness of smooth solutions
to the Euler/Allen-Cahn system. In section 4, we prove that the
solution of the Navier-Stokes/Allen-Cahn equations converges to that
of the Euler/Allen-Cahn equations under the condition \eqref{Kato}.
In the last section, we show some results for related systems.

\section{Preliminaries and Main Results}

Let $C^{\infty}_{0,div}(\Omega)$ be the space of all divergence free
vectors in $(C^{\infty}_{0}(\Omega))^n$ ($n=2,3$). We denote by $H$
the closure of $C^{\infty}_{0,div}(\Omega)$ in $(L^2(\Omega))^n$.
Moreover, we set \beee
 && V=(H^1_0(\Omega))^n\cap H,  \hspace{3cm}
  V_2=(H^2(\Omega))^n\cap V,\\
&&\mathcal{W}=\left\{\phi\in C^{\infty}(\Omega),\
 \partial_\textbf{n}\phi|_{\partial\Omega}=0 \right\},\quad \quad
\Phi_s=\mbox{the closure of}\, \mathcal{W} \,\mbox{in} \,
H^s(\Omega),\, s\in \mathbb{N}^+,\,\\
&&\mathcal {O}=\left\{ {\bf v}\in (C^{\infty}(\Omega))^n,
\nabla\cdot{\bf v}=0\,\,\mbox{in}\ \Omega, \, {\bf v}\cdot{\bf
n}=0\,\,\mbox{on}\ \partial\Omega  \right\},\\
&& X_0=\mbox{the closure of}\, \mathcal{O} \,\mbox{in} \,
(L^2(\Omega))^n,\quad \quad X_s = X_0\cap(H^s(\Omega))^n, s\in
\mathbb{N}^+ .\eeee Obviously, $H\subset X_0$ (continuous
embedding), $X_s\subset\subset X_0$ (compact embedding) for $s\geq
1$.

Let $\Gamma_\delta$ be the boundary strip of width $\delta$. We
shall take $\delta=c\nu$ with $c>0$ being a small constant. We
denote $|\cdot|$ and $(\cdot,\cdot)$ as the norm and scalar product
in $L^2(\Omega)$ or $(L^2(\Omega))^n$. For any positive integer $s$,
we take $|\cdot|_s$ and $((\cdot,\cdot))_s$ as the norm and scalar
product in $H^s(\Omega)$ or $(H^s(\Omega))^n$,
$((f,g))_s=\sum_{|\alpha|\leq s}\left(D^\alpha f,D^\alpha g\right)$
where $D^\alpha$ is a multi-index derivation. We also denote
$|\cdot|_X$ as the norm in other Banach spaces $X$.

Following almost the same arguments as in \cite{LL95}, for any fixed
$\nu>0$, using the important energy law \eqref{basic energy law}, we
have the following well-posedness result for the
Navier-Stokes/Allen-Cahn system:
\bt\label{thm1} Assume that $n=2,3$. ${\bf u}_0\in H$ and $\phi_0\in
\Phi_1$. Then the system \eqref{Navier-Stokes}-\eqref{BC1} has a
global weak solution $({\bf u},\phi)$ such that for all
$T\in(0,+\infty)$ \beee &&{\bf u}\in
L^\infty(0,T;H)\cap L^2(0,T;V), \\
&& \phi\in L^\infty(0,T;\Phi_1)\cap L^2(0,T;\Phi_2).\eeee Moreover,
if ${\bf u}_0\in V$ and $\phi_0\in \Phi_2$, then in the $2$D case,
system \eqref{Navier-Stokes}-\eqref{BC1} admits a unique global
strong solution such that for all $T\in(0,+\infty)$
 \bee &&{\bf u}\in
L^\infty(0,T;V)\cap L^2(0,T;V_2),\label{strong1}
\\ && \phi\in L^\infty(0,T;\Phi_2)\cap
L^2(0,T;\Phi_3); \label{strong-soln}\eee while in the $3$D case,
there exists a $T_\nu>0$ depending on $|{\bf u}_0|_1$, $|\phi_0|_2$
and $\nu$ such that \eqref{Navier-Stokes}-\eqref{BC1} has a unique
strong solution in $[0,T_\nu)$ and
\eqref{strong1}-\eqref{strong-soln} hold for $T<T_\nu$. \et

For weak solution $({\bf u},\phi)$ to system
\eqref{Navier-Stokes}-\eqref{BC1}, we easily deduce from energy law
\eqref{basic energy law} that for any $T>0$, \be
\nu\int^T_0|\nabla{\bf u}|^2(t)dt\leq C, \quad\quad
\int^T_0|\phi|^2_2(t)dt\leq C, \label{bound}\ee where $C$ depends on
${\bf u}_0, \phi_0$ and $|\Omega|$.

We also prove the local existence of smooth solutions to
Euler/Allen-Cahn system.
 \bt\label{thm3} Assume that $n=2,3$. ${\bf v}_0\in X_s$, $\psi_0\in \Phi_{s+1}$, where $s\geq
 3$ is an integer. Then there exists a small constant $T_*>0$
 depending on $|{\bf v}_0|_s$ and $|\psi_0|_{s+1}$, such that
the system \eqref{Euler}-\eqref{BC2} admits a smooth solution $({\bf
v}, P_2, \psi)$ on $[0,T_*]$ satisfying \be {\bf v}\in
L^\infty(0,T_*;X_s),\quad P_2\in L^\infty(0,T_*; H^{s+1}(\Omega)),
\quad \psi\in L^\infty(0,T_*;\Phi_{s+1}). \ee Moreover, $({\bf v},
\psi)$ is unique in the sense of
$\left(L^\infty(0,T_*;L^2(\Omega)),L^\infty(0,T_*;H^1(\Omega))\right)$-norm.
 \et Finally, we show the convergence
of solutions of Navier-Stokes/Allen-Cahn system to that of
Euler/Allen-Cahn system.
\bt\label{thm2} Assume that $n=2,3$. Let $({\bf u},\phi)$ be a
global weak solution to \eqref{Navier-Stokes}-\eqref{BC1}, and
$({\bf v},\psi)$ a local smooth solution to
\eqref{Euler}-\eqref{BC2} on $\Omega\times[0,T_*]$. Assume that \be
\nu\int^{T_*}_0|\nabla{\bf
u}|^2_{L^2(\Gamma_{c\nu})}(t)dt\rightarrow 0 \ \mbox{as}\
\nu\rightarrow 0.\label{condition2}\ee Then \be|{\bf u}(t)-{\bf
v}(t)|^2+\lambda|\phi(t)-\psi(t)|_{1}^2\rightarrow 0 \ \mbox{as}\
\nu\rightarrow 0, \quad \mbox{for a.e.}\ t\in[0,T_*].\label{con}\ee
Moreover, if $n=2$ and $({\bf u},\phi)$ is a global strong solution
to \eqref{Navier-Stokes}-\eqref{BC1}, then the assumption
\eqref{condition2} is automatically satisfied, and the convergence
rate in \eqref{con} is $c\nu^{1/2}$.\et
\br In general, we do not know whether the weak solution $({\bf
u},\phi)$ satisfies the assumption \eqref{condition2}. We only know
$\nu\int^{T_*}_0|\nabla{\bf u}|^2_{L^2(\Gamma_{c\nu})}(t)dt$ is
bounded for weak solutions (see \eqref{bound}). But in 2D, the
global strong solutions exist. If $({\bf u},\phi)$ is a global
strong solutions, i.e. $\|\nabla{\bf
u}\|_{L^\infty(0,T;L^2(\Omega))}$ is bounded and its upper bound is
independent of $\nu$, then ${\bf u}$ satisfies \eqref{condition2}
obviously. \er \br The above theorems remain true if we replace the
Neumann boundary condition for $\phi$ and $\psi$ with homogeneous
Dirichlet boundary condition. \er

\section{Local smooth solutions to Euler/Allen-Cahn system }
In this section, we give the proof of Theorem \ref{thm3} in the 3D
case. The 2D case can be proved similarly. We shall apply a modified
Galerkin method, following the spirit of \cite{LL95}. We choose an
orthonormal complete basis $\{{\bf w}_k\}^\infty_{k=1}\subset X_s$
$(s\geq 3)$ used in \cite{Temam1975} for Euler equation, which
satisfies
%Obviously, $X_s\subset\subset X_0$ and $X_s$ is dense in $X_0$. By
%Lax-Milgram theorem, for any ${\bf g}\in X_0$, there exists a unique
%${\bf w({\bf g})}\in X_s$ such that \[(({\bf w({\bf g})},{\bf
%v}))_s=({\bf g},{\bf v}),\ \ \forall {\bf v}\in X_s.\] It is easy to
%verify that the linear mapping ${\bf g}\rightarrow {\bf w({\bf g})}$
%is a compact self-adjoint operator in $X_0$ and has a complete
%family of eigenvectors ${\bf w}_k\in X_s$ such that
%
\bee \left\{
\begin{array}{c}
(({\bf w}_k,{\bf
v}))_s=\lambda_k({\bf w}_k,{\bf v}),\ \ \forall {\bf v}\in X_s.\\
({\bf w}_k,{\bf w}_j)=\delta_{kj}, i.e., (({\bf w}_k,{\bf
w}_j))_s=\lambda_k\delta_{kj}.\end{array}\right.\eee Here,
$0<\lambda_1\leq\lambda_2\leq\cdots\leq\lambda_n\leq\cdots$ and
$\lambda_n\rightarrow\infty$ as $n\rightarrow\infty$.

For fixed integer $m>0$, denote $X_{sm}$ as the space spanned by
$\{{\bf w}_1,{\bf w}_2,\cdots,{\bf w}_m\}$. Let $P_m:X_0\rightarrow
X_{sm}$ be the orthonormal projection. We consider the approximate
problem: \bee &&\frac{\partial{\bf v}_m}{\partial t}=P_m\left(-({\bf
v}_m\cdot\nabla){\bf
v}_m-\lambda\nabla\cdot\left(\nabla\psi_m\otimes\nabla\psi_m\right)\right),\quad \ {\bf v}_m\in X_{sm},\label{approx-1}\\
&&\frac{\partial\psi_m}{\partial t}+({\bf
v}_m\cdot\nabla)\psi_m=\gamma\left(\Delta\psi_m-f(\psi_m)\right), \label{approx-2}\\
&&{\bf v}_m(0,x)=P_m{\bf v}_0(x),\quad \psi_m(0,x)=\psi_0(x),\\
&&{\bf v}_m\cdot{\bf n}|_{\partial\Omega}=0,\quad
\partial_\textbf{n}\psi_m|_{\partial\Omega}=0.
\label{approx-3}\eee Here, we only project \eqref{Euler} into finite
dimensional space $X_{sm}$, but do not project \eqref{psi}. We shall
see later in \eqref{10} that, the approximate solutions $({\bf
v}_m,\psi_m)$ still satisfies energy law. This guarantees the global
existence of $({\bf v}_m,\psi_m)$.

We apply Schauder fixed point theorem to deduce the local existence
of solutions to \eqref{approx-1}-\eqref{approx-3}. The arguments are
similar to that in \cite[Therem 2.1]{LL95}. For readers'
convenience, we give the outline of the proof here. Consider a
closed, convex subset in $(C[0,T_0])^m$: {\small
\[D=\left\{(g^1_m,g^2_m,\cdots,g^m_m)\left|
\left(\sum^{m}_{i=1}|g^i_m(t)|^2\right)^{1/2}\leq M,\, 0\leq t\leq
T_0;\, g^i_m(0)=({\bf v}_0,{\bf w}_i),
i=1,2,\cdots,m.\right.\right\}\]} where $M$ and $T_0$ are positive
constants to be determined later. Given a
$(g^1_m,g^2_m,\cdots,g^m_m)\in D$, we get a ${\bf
v}_m(x,t)=\sum^{m}_{i=1}g^i_m(t){\bf w}_i(x)$ satisfying $|{\bf
v}_m|_{L^\infty(\Omega)}\leq c_mM$ for $t\in[0,T_0]$, where $c_m$ is
a constant depending on $m$. With this ${\bf v}_m$ in
\eqref{approx-2}, there exists a regular solution $\psi_m$ on
$\Omega\times[0,T_0]$, and \ube|\nabla\psi_m|^2(t)\leq
e^{c_m^2M^2t}\left(|\nabla\psi_0|^2+2\int_\Omega
F(\psi_0)dx\right)=:e^{c_m^2M^2t}C(|\psi_0|_1,|\Omega|),\quad
t\in[0,T_0].\uee
Then we substitute this $\psi_m$ into \eqref{approx-1} and look for
a solution ${\bf \tilde{v}}_m=\sum^{m}_{i=1}\tilde{g}^i_m(t){\bf
w}_i(x)$. It is standard that \eqref{approx-1} is equivalent to an
ODE system of $\tilde{g}^i_m(t)$ $(i=1,2,\cdots,m)$, with initial
data $\tilde{g}^i_m(0)=({\bf v}_0,{\bf w}_i)$. There exists a local
solution $(\tilde{g}^1_m,\tilde{g}^2_m,\cdots,\tilde{g}^m_m)\in
(C^1[0,T_0])^m$. Moreover, it turns out that
{\small\bee\left(\sum^{m}_{i=1}|\tilde{g}^i_m(t)|^2\right)^{1/2}&\leq&\left(\sum^{m}_{i=1}|\tilde{g}^i_m(0)|^2\right)^{1/2}
+\frac{\lambda
C(m,|\psi_0|_1,\Omega)}{M^2}\left(e^{c_m^2M^2t}-1\right)\non\\
&=&\left(\sum^{m}_{i=1}|({\bf v}_0,{\bf w}_i)|^2\right)^{1/2}
+\frac{\lambda
C(m,|\psi_0|_1,\Omega)}{M^2}\left(e^{c_m^2M^2t}-1\right)
\label{g}\eee} for $t\in[0,T_0]$. Let $M=2\left(\sum^{m}_{i=1}|({\bf
v}_0,{\bf w}_i)|^2\right)^{1/2}+2$. Then there exits a small $T_0$
such that the right hand side of \eqref{g} is not more than $M$.
Therefore,
${\small\left(\sum^{m}_{i=1}|\tilde{g}^i_m(t)|^2\right)^{1/2}\leq
M}$ for $t\in[0,T_0]$. Thus
$(\tilde{g}^1_m,\tilde{g}^2_m,\cdots,\tilde{g}^m_m)\in D$. And the
mapping
$\mathcal{L}:(g^1_m,g^2_m,\cdots,g^m_m)\mapsto(\tilde{g}^1_m,\tilde{g}^2_m,\cdots,\tilde{g}^m_m)$
is a compact operator in $(C[0,T_0])^m$. Applying Schauder fixed
point theorem, $\mathcal{L}$ has a fixed point. Thus we obtain the
following lemma: \bl\label{lemma2} There exists a $T_0>0$, depending
on ${\bf v}_0, \psi_0, m$ and $\Omega$, such that
\eqref{approx-1}-\eqref{approx-3} has a weak solution $({\bf
v}_m,\psi_m)$ in $\Omega\times[0,T_0]$.\el
In order to show the global existence of weak solution to
\eqref{approx-1}-\eqref{approx-3}, we establish a priori estimates
next. Assume $({\bf v}_m,\psi_m)$ is a weak solution to
\eqref{approx-1}-\eqref{approx-3} in $\Omega\times[0,T]$ for certain
$T>0$. Taking the inner product of \eqref{approx-1} with ${\bf v}_m$
in $L^2(\Omega)$, using the fact that $\left(({\bf
v}_m\cdot\nabla){\bf v}_m, {\bf v}_m\right)=0$ and
$\nabla\cdot\left(\nabla\psi_m\otimes\nabla\psi_m\right)=\nabla\left(\frac{|\nabla\psi_m|^2}{2}\right)+
\Delta\psi_m\nabla\psi_m$, we get \be\frac{1}{2}\frac{d}{dt}|{\bf
v}_m|^2=-\lambda\left(\Delta\psi_m\nabla\psi_m, {\bf v}_m\right).
\ee Taking the inner product of \eqref{approx-2} with
$\gamma\left(\Delta\psi_m-f(\psi_m)\right)$ in $L^2(\Omega)$, we
have \be
\frac{d}{dt}\left(\frac{\lambda}{2}|\nabla\psi_m|^2+\lambda\int_\Omega
F(\psi_m)dx\right)+\lambda\gamma|\Delta\psi_m-f(\psi_m)|^2=\lambda\left(({\bf
v}_m\cdot\nabla)\psi_m,\Delta\psi_m\right),\ee where we used
$\left(({\bf v}_m\cdot\nabla)\psi_m,f(\psi_m)\right)=\int_\Omega{\bf
v}_m\cdot\nabla F(\psi_m)=0$ where
$F(\psi_m)=\frac{1}{4\varepsilon^2}(\psi_m^2-1)^2$. Adding the above
two resultant, we obtain \be\frac{d}{dt}\left(\frac{1}{2}|{\bf
v}_m|^2+\frac{\lambda}{2}|\nabla\psi_m|^2+\lambda\int_\Omega
F(\psi_m)dx\right)+\lambda\gamma|\Delta\psi_m-f(\psi_m)|^2=0.\label{10}\ee
Hence, \be \sup_{0\leq t\leq T}\left(|{\bf
v}_m|^2(t)+\lambda|\nabla\psi_m|^2(t)+\lambda\int_\Omega
F(\psi_m)dx\right) \leq |{\bf
v}_0|^2+\lambda|\nabla\psi_0|^2+\lambda\int_\Omega F(\psi_0)dx. \ee
By Lemma \ref{lemma2} and the above inequality, we have \bt For any
integer $m>0$, ${\bf v}_0\in X_0$, $\psi_0\in \Phi_1$, system
\eqref{approx-1}-\eqref{approx-3} has a weak solution $({\bf v}_m,
\psi_m)$ in $\Omega\times[0,\infty)$. \et

In what follows, we do higher order energy estimates, so that we can
pass to limits as $m\rightarrow \infty$. Firstly, since ${\bf
v}_m\in X_{sm}$, we can rewrite \eqref{approx-1} as
\be\frac{\partial{\bf v}_m}{\partial t}=-({\bf v}_m\cdot\nabla){\bf
v}_m-\lambda\Delta\psi_m\nabla\psi_m-\nabla q_m,\label{vm-eqn}\ee
for a certain $q_m$. By taking divergence operator on both sides of
\eqref{vm-eqn} on $\Omega$, and taking scalar product of both sides
of \eqref{vm-eqn} with ${\bf n}$ on $\partial\Omega$, we know that
$q_m$ satisfies \beee \Delta q_m
&=&-\lambda\nabla\cdot(\Delta\psi_m\nabla\psi_m)-\sum_{i,j}\nabla_j{\bf
v}_{mi}\nabla_i{\bf v}_{mj},\ \mbox{in}\ \Omega \\
\frac{\partial q_m}{\partial{\bf
n}}&=&-\lambda\Delta\psi_m\frac{\partial\psi_m}{\partial{\bf
n}}+\sum_{i,j}{\bf v}_{mi}{\bf v}_{mj}\Gamma_{ij}, \ \mbox{on}\
\partial\Omega\eeee
where the function $\Gamma_{ij}$ depends only on $\partial\Omega$
(see \cite[Lemma 1.1]{Temam1975}).
 By elliptic estimates (see also
\cite[Lemma 1.2]{Temam1975}), for $s\geq \frac{5}{2}$,
\be\label{q-m}|\nabla q_m|_s\leq
C\left(|\Delta\psi_m|_s|\nabla\psi_m|_s+|{\bf v}_m|^2_s\right).\ee

Taking inner product of \eqref{vm-eqn} with $\lambda_kg_{km}{\bf
w}_k$ in $L^2(\Omega)$, and adding in $k$, $k=1,2\cdots,m$, we
obtain \be\label{6} \frac{1}{2}\frac{d}{dt}|{\bf v}_m|^2_s=((-({\bf
v}_m\cdot\nabla){\bf v}_m-\lambda\Delta\psi_m\nabla\psi_m-\nabla
q_m, {\bf v}_m ))_s,\ee

For $s\geq 3$, we estimate the righthand side of \eqref{6} to get
\bee \frac{1}{2}\frac{d}{dt}|{\bf v}_m|^2_s&\leq&\left|(({\bf
v}_m\cdot\nabla{\bf v}_m, {\bf
v}_m))_s\right|+\left|\lambda((\Delta\psi_m\nabla\psi_m,{\bf
v}_m))_s\right|+\left|((\nabla q_m,{\bf v}_m))_s\right|\non\\
&\leq&C|{\bf v}_m|^3_s+C\lambda|\Delta\psi_m|_s|\nabla\psi_m|_s|{\bf
v}_m|_s+|\nabla q_m|_s|{\bf v}_m|_s\non\\
&\leq&C|{\bf
v}_m|^3_s+\frac{\lambda\gamma}{8}|\Delta\psi_m|_s^2+C|\nabla\psi_m|_s^4+C|{\bf
v}_m|_s^4+C\left(|\Delta\psi_m|_s|\nabla\psi_m|_s+|{\bf v}_m|_s^2
\right)|{\bf v}_m|_s      \non\\
&\leq&\frac{\lambda\gamma}{4}|\Delta\psi_m|_s^2+C|\nabla\psi_m|_s^4+C|{\bf
v}_m|_s^4+C, \label{um}\eee where we used $\left(({\bf
v}_m\cdot\nabla)\nabla^s{\bf v}_m, \nabla^s{\bf v}_m\right)=0$ and
\eqref{q-m}. \br In the 2D case, $H^{1+\varepsilon}(\Omega)\subset
L^\infty(\Omega),\, \forall \varepsilon>0$. $\left|(({\bf
v}_m\cdot\nabla{\bf v}_m, {\bf v}_m))_s\right|$ in \eqref{um} can be
controlled by $|{\bf v}_m|_s^\alpha$ optimally, where $2<\alpha<3$.
This term makes it difficult to obtain global existence of smooth
solutions.\er

Applying $\nabla^i$ ($i=1,2,\cdots, s;\, s\geq 3$) to
\eqref{approx-2} (here, $\nabla^i=\nabla(-\Delta)^{(i-1)/2}$ if $i$
is odd, $\nabla^i=(-\Delta)^{i/2}$ if $i$ is even), then taking the
$L^2$ inner product of the resulting equation with
$\lambda\nabla^i\Delta\psi_m$, and using the boundary condition
$\frac{\partial\psi_m}{\partial{\bf n}}|_{\partial\Omega}=0$ and
higher order natural boundary conditions such as
$\frac{\partial\Delta\psi_m}{\partial{\bf n}}|_{\partial\Omega}=0$,
we obtain \ube
\frac{\lambda}{2}\frac{d}{dt}|\nabla^{i+1}\psi_m|^2+\lambda\gamma|\nabla^i\Delta\psi_m|^2
=\lambda\left(\nabla^i({\bf
v}_m\cdot\nabla)\psi_m,\nabla^i\Delta\psi_m\right)+\lambda\gamma\left(\nabla^if(\psi_m),\nabla^i\Delta\psi_m\right).\uee
Adding the above resultant in $i$ $(i=1,2,\cdots,s;\,s\geq 3)$, we
get\be
\frac{\lambda}{2}\frac{d}{dt}|\nabla\psi_m|_s^2+\lambda\gamma|\Delta\psi_m|_s^2=\lambda\sum^s_{i=1}\left(\nabla^i({\bf
v}_m\cdot\nabla)\psi_m,\nabla^i\Delta\psi_m\right)+\lambda\gamma\sum^s_{i=1}\left(\nabla^if(\psi_m),\nabla^i\Delta\psi_m\right).\label{psi-s}\ee
 Next we estimate the righthand side terms. \bee&&
\lambda\left|\sum^s_{i=1}\left(\nabla^i({\bf
v}_m\cdot\nabla)\psi_m,\nabla^i\Delta\psi_m\right)\right|\non\\
&\leq& \sum^s_{i=1}\lambda\left(|\nabla^i{\bf
v}_m||\nabla\psi_m|_{L^\infty}+|{\bf
v}_m|_{L^\infty}|\nabla^{i+1}\psi_m|\right)|\nabla^i\Delta\psi_m|\non\\
&\leq& \sum^s_{i=1}\lambda\left(|\nabla^i{\bf
v}_m||\nabla\psi_m|_2+|{\bf
v}_m|_2|\nabla^{i+1}\psi_m|\right)|\nabla^i\Delta\psi_m|\non \\
&\leq&\sum^s_{i=1}\left[\frac{\lambda\gamma}{4}|\nabla^i\Delta\psi_m|^2+C|\nabla^i{\bf
v}_m|^4+C|\nabla\psi_m|_2^4+C|{\bf
v}_m|_2^4+C|\nabla^{i+1}\psi_m|^4 \right]\non\\
&\leq&\frac{\lambda\gamma}{8}|\Delta\psi_m|_s^2+C|{\bf
v}_m|_s^4+C|\nabla\psi_m|_s^4, \eee
\bee
&&\lambda\gamma\left(\nabla^if(\psi_m),\nabla^i\Delta\psi_m\right)\non\\
&=&\frac{4\lambda\gamma}{\varepsilon^2}\left(\nabla^i(\psi_m^3),\nabla^i\Delta\psi_m\right)
-\frac{4\lambda\gamma}{\varepsilon^2}\left(\nabla^i\psi_m,\nabla^i\Delta\psi_m\right)\non\\
&=&\frac{12\lambda\gamma}{\varepsilon^2}\left(\psi_m^2\nabla^i\psi_m,\nabla^i\Delta\psi_m\right)
+\sum\left(C\nabla^{i_1}\psi_m\nabla^{i_2}\psi_m\nabla^{i_3}\psi_m,\nabla^i\Delta\psi_m
\right)\non\\
&&-\frac{4\lambda\gamma}{\varepsilon^2}\left(\nabla^i\psi_m,
\nabla^i\Delta\psi_m\right)\non\\
&=:&I_1+I_2+I_3, \eee
where $\Sigma$ represents the sum with respect to integers $i_1, i_2, i_3$ satisfying $0\leq i_1\leq i_2\leq i_3\leq i-1$ and
$i_1+i_2+i_3=i$.

Using the fact that $\psi_m$ is uniformly bounded in
$L^\infty(0,T;H^1(\Omega))$ and $H^1\subset L^6$, we get\bee
|I_1|&\leq&C|\psi_m|_{L^6}^2|\nabla^i\psi_m|_{L^6}|\nabla^i\Delta\psi_m|\non\\
&\leq&C|\psi_m|_1^2|\nabla^i\psi_m|_1|\nabla^i\Delta\psi_m|\non\\
&\leq&C|\nabla\psi_m|_i|\nabla^i\Delta\psi_m|\non\\
&\leq&\frac{\lambda\gamma}{8}|\nabla^i\Delta\psi_m|^2+C|\nabla\psi_m|_i^2,
\eee
\bee |I_2|&\leq&C\sum|\nabla^{i_1}\psi_m|_{L^6}|\nabla^{i_2}\psi_m|_{L^6}|\nabla^{i_3}\psi_m|_{L^6}
|\nabla^i\Delta\psi_m|\non\\
&\leq&C\sum|\nabla^i\Delta\psi_m|\prod^{3}_{k=1}\left(C|\nabla\psi_m|^{1-a_k}|\nabla^i\Delta\psi_m|^{a_k}+C|\nabla\psi_m|\right)\non\\
&\leq&C\sum|\nabla^i\Delta\psi_m|\left(|\nabla^i\Delta\psi_m|^{\frac{i}{i+1}}+1\right),\non\\
&\leq&\frac{\lambda\gamma}{8}|\nabla^i\Delta\psi_m|^2+C, \eee where
we used the Gagliardo--Nirenberg interpolation inequalities (cf.
\cite{Z04}): \ube |\nabla^{i_k}\psi_m|_{L^6}\leq
C|\nabla\psi_m|^{1-a_k}|\nabla^i\Delta\psi_m|^{a_k}+C|\nabla\psi_m|,
\quad a_k=\frac{i_k}{i+1},\quad k=1,2,3.\uee
%|\nabla^{i_2}\psi|_{L^4}&\leq&C|\nabla\psi|^{1-a_2}|\nabla^i\Delta\psi|^{a_2},
%\quad a_2=\frac{3}{i+1}\left(\frac{1}{4}+\frac{i_2-1}{3}\right),\\
%|\nabla^{i_3}\psi|_{L^4}&\leq&C|\nabla\psi|^{1-a_3}|\nabla^i\Delta\psi|^{a_3},
%\quad a_3=\frac{3}{i+1}\left(\frac{1}{4}+\frac{i_3-1}{3}\right).
%\eeee
and $a_1+a_2+a_3=\frac{i}{i+1}<1$.
\be |I_3|\leq C|\nabla^i\psi_m||\nabla^i\Delta\psi_m|\leq
\frac{\lambda\gamma}{8}|\nabla^i\Delta\psi_m|^2+C|\nabla^i\psi_m|^2.\ee
Using the above estimates, we conclude from \eqref{psi-s} that \be
\frac{\lambda}{2}\frac{d}{dt}|\nabla\psi_m|^2_s+\frac{\lambda\gamma}{2}|\Delta\psi_m|_s^2
\leq C|\nabla\psi_m|_{s}^4+C|{\bf v}_m|_s^4+C. \label{high-psi} \ee
Adding \eqref{um} and \eqref{high-psi}, we conclude \be\label{9}
\frac{1}{2}\frac{d}{dt}\left(|{\bf
v}_m|_s^2+\lambda|\nabla\psi_m|^2_s\right)+\frac{\lambda\gamma}{4}|\Delta\psi_m|_s^2
\leq C|{\bf v}_m|_s^4+C|\nabla\psi_m|_s^4+C. \ee
Denote $Y_m(t)=|{\bf v}_m|_s^2+\lambda|\nabla\psi_m|_s^2$. \eqref{9}
can be rewritten as
\[\frac{d}{dt}Y_m(t)\leq C_1Y_m^2(t)+C_2,\]
\[Y_m(0)=|{\bf v}_m|_s^2(0)+\lambda|\nabla\psi_m|^2_s(0)\leq |{\bf v}_0|_s^2+\lambda|\nabla\psi_0|^2_s.\]
Hence there exists a $T_*$ depending only on $C_1, C_2,$ and $|{\bf
v}_0|_s$, $|\psi_0|_{s+1}$, such that $Y_m(t)\leq N$ on $[0,T_*]$,
where $N>0$ is a constant independent of $m$. Therefore, as
$m\rightarrow +\infty$, \bee &&{\bf v}_m \ \mbox{remains bounded in
}\,
L^\infty(0,T_*;X_s),\label{bdd-1}\\
&&\psi_m \  \mbox{remains bounded in }\,
L^\infty(0,T_*;\Phi_{s+1})).\label{bdd-2}\eee
By the equation \eqref{approx-1} and \eqref{approx-2}, it is easy to
know that \bee&&\frac{\partial{\bf v}_m}{\partial t} \ \mbox{remains
bounded in }\, L^\infty(0,T_*;X_{s-1}),\label{bdd-3}\\
&&\frac{\partial\psi_m}{\partial t} \  \mbox{remains bounded in }\,
L^\infty(0,T_*;\Phi_{s-1}).\label{bdd-4}\eee
Using \eqref{bdd-1}-\eqref{bdd-4} and the standard compact embedding
theorem, the passage to the limit is standard. We obtain the
existence of ${\bf v}\in L^\infty(0,T_*;X_s)$, $\psi\in
L^\infty(0,T_*;\Phi_{s+1})$ such that \beee
&&\left(\frac{\partial{\bf v}}{\partial t},
{\bm\xi}\right)=\left(-({\bf v}\cdot\nabla){\bf
v}-\lambda\nabla\cdot\left(\nabla\psi\otimes\nabla\psi\right),{\bm\xi}\right),\quad \forall {\bm\xi}\in X_0, \quad \mbox{a.e.}\ t\in[0,T_*],\label{5}\\
&&(\frac{\partial\psi}{\partial t},\varphi)+\left(({\bf
v}_m\cdot\nabla)\psi,\varphi\right)=-\gamma\left(\nabla\psi,\nabla\varphi\right)-\left(f(\psi),\varphi\right),
\quad \forall \varphi\in H^1(\Omega),\quad \mbox{a.e.}\ t\in[0,T_*].
\eeee The proof of global existence of smooth solutions to
\eqref{Euler}-\eqref{BC2} is finished.

It is standard to prove the uniqueness of smooth solutions to
\eqref{Euler}-\eqref{BC2}. We omit it here. The proof of Theorem \ref{thm3} is complete.

\section{Vanishing viscosity limit}

In this section, we prove Theorem \ref{thm2}. In the proof, $C$
is a constant independent of $\nu$.

To deal with the mismatch between the boundary condition of ${\bf
u}$ and that of ${\bf v}$, we use the boundary layer function
$\bm\theta$ constructed in \cite{Kato1984}:
$\bm\theta=div\left(\chi\left(\frac{\rho}{\delta}\right){\bf
A}\right)$ where $\rho=dist({\bf x},\partial\Omega)$ is a distance
function, $\delta$ is a small positive constant, and $\chi$ is a
smooth cut-off function
$\chi:\mathbb{R}^{+}\rightarrow\mathbb{R}^{+}$ such that
\[\chi(0)=1,\quad \chi(r)=0 \ \mbox{for}\ r\geq 1.\]
${\bf A}$ is a skew-symmetric matrix defined on $\Omega\times[0,T]$,
such that \be div{\bf A}={\bf v} \ \mbox{on} \
\partial\Omega,\ {\bf A}=0 \ \mbox{on}\ \partial\Omega.\ee
Obviously, $\bm\theta$ is supported on $\Gamma_\delta$. It was
proved in \cite{Kato1984} that $\bm\theta$ satisfies
\[div \bm\theta=0\ \mbox{in}\ \Omega, \quad \bm\theta={\bf v}\ \mbox{on}\ \partial\Omega,\]
and $\bm\theta$ has the following estimates: \bee\label{theta-1}
&&|\bm\theta|_{L^\infty}\leq C, \quad |\bm\theta|\leq
C\delta^{\frac{1}{2}},\quad |\bm\theta|_{L^4}\leq
C\delta^{\frac{1}{4}},
\quad |\bm\theta_t|\leq C\delta^{\frac{1}{2}},\\
&&|\nabla\bm\theta|_{L^\infty}\leq C\delta^{-1},\quad
|\nabla\bm\theta|\leq C\delta^{-\frac{1}{2}},\quad
|\rho\nabla\bm\theta|_{L^\infty}\leq C, \\
&&|\rho^2\nabla\bm\theta|_{L^\infty}\leq C\delta,\quad
|\rho\nabla\bm\theta|\leq C\delta^{\frac{1}{2}},\label{theta-2} \eee
provided that ${\bf v}$ has the regularity: \be |\nabla^2{\bf
v}|_{L^\infty}\leq C,\quad |\nabla{\bf v}_t|_{L^\infty}\leq C,\quad
|{\bf v}_t|_{L^\infty}\leq C. \ee

Suppose $({\bf u},\phi)$ is a global weak solution to
\eqref{Navier-Stokes}-\eqref{BC1}, and $({\bf v},\psi)$ is a local
smooth solution to \eqref{Euler}-\eqref{BC2} on
$\Omega\times[0,T_*]$.
 Denote ${\bf w}={\bf u}-{\bf v}+\bm\theta$, $\zeta=\phi-\psi$, then
${\bf w}, \zeta$ satisfy the following equations: \bee {\bf
w}_t+({\bf u}\cdot\nabla){\bf w}+\nabla q &=& \nu\Delta{\bf u}-({\bf
w}\cdot\nabla){\bf v}+(\bm\theta\cdot\nabla){\bf v}+({\bf
u}\cdot\nabla)\bm\theta+\bm\theta_t\non\\
&&-\lambda\Delta\zeta\nabla\phi-\lambda\Delta\psi\nabla\zeta,\label{w-eqn}\\
\nabla\cdot {\bf w} &=& 0,\\
\zeta_t+({\bf v}\cdot\nabla)\zeta+({\bf
w}\cdot\nabla)\phi &=& \gamma\Delta\zeta+(\bm\theta\cdot\nabla)\phi-\frac{4\gamma}{\varepsilon^2}\zeta(\phi^2+\phi\psi+\psi^2)+\frac{4\gamma}{\varepsilon^2}\zeta,\label{zeta-eqn}\\
{\bf w}|_{t=0}=\bm\theta(x,0),&&
\zeta|_{t=0}=0,\ \ x\in\Omega,\\
{\bf w}|_{\partial\Omega}(x,t)=0,&&
\partial_\textbf{n}\zeta|_{\partial\Omega}(x,t)=0,\ \ t\geq 0,
 \eee
where
$q=P_1-P_2+\frac{\lambda}{2}|\nabla\phi|^2-\frac{\lambda}{2}|\nabla\psi|^2$.

Our goal is to estimate $|{\bf w}|(t)$ and $|\zeta|_{H^1}(t)$, we
hope they go to zero, equivalently,
$|{\bf u}-{\bf v}|(t)$ and $|\phi-\psi|_{H^1}(t)$ converge to zero
as $\nu$ tends to zero.

Multiplying \eqref{w-eqn} by ${\bf w}$ and integrating over
$\Omega$, we obtain \bee \frac{1}{2}\frac{d}{dt}|{\bf
w}|^2(t)+\nu|\nabla{\bf u}|^2&=&\nu\left(\nabla{\bf u},\nabla({\bf
v}-\bm\theta)\right)-\left(({\bf w}\cdot\nabla){\bf v},{\bf
w}\right)+\left((\bm\theta\cdot\nabla){\bf v},{\bf w}\right)+\left(({\bf u}\cdot\nabla)\bm\theta,{\bf u}\right)\non\\
&&-\left(({\bf u}\cdot\nabla)\bm\theta,{\bf
v}\right)+(\bm\theta_t,{\bf
w})-\lambda\left(\Delta\zeta\nabla\phi,{\bf
w}\right)-\lambda\left(\Delta\psi\nabla\zeta,{\bf w}\right).
\label{estimate-w} \eee
Multiplying \eqref{zeta-eqn} by $\lambda\zeta$ and integrating over
$\Omega$, we get \bee &&\frac{\lambda}{2}\frac{d}{dt}|\zeta|^2+
\lambda\gamma|\nabla\zeta|^2+\frac{4\lambda\gamma}{\varepsilon^2}\int_\Omega\zeta^2\left(\phi^2+\phi\psi+\psi^2\right)dx\non\\
&=&\frac{4\lambda\gamma}{\varepsilon^2}|\zeta|^2-\lambda({\bf
w}\cdot\nabla\phi,\zeta)+\lambda({\bm\theta\cdot\nabla}\phi,\zeta).
\label{estimate-zeta1}\eee Notice that
$\phi^2+\phi\psi+\psi^2=\frac{1}{2}\left(\phi^2+\psi^2+(\phi+\psi)^2\right)\geq
0$, thus
$\int_\Omega\zeta^2\left(\phi^2+\phi\psi+\psi^2\right)dx\geq 0$.
Multiplying \eqref{zeta-eqn} by $\lambda\Delta\zeta$ and integrating
over $\Omega$, we obtain \bee
\frac{\lambda}{2}\frac{d}{dt}|\nabla\zeta|^2+\lambda\gamma|\Delta\zeta|^2&=&\lambda\left(({\bf
v}\cdot\nabla)\zeta,\Delta\zeta\right)+\lambda\left(({\bf
w}\cdot\nabla)\phi,\Delta\zeta\right)-\lambda\left((\bm\theta\cdot\nabla)\phi,\Delta\zeta\right)\non\\
&&+\frac{4\lambda\gamma}{\varepsilon^2}\left(\zeta(\phi^2+\phi\psi+\psi^2),\Delta\zeta\right)+\frac{4\lambda\gamma}{\varepsilon^2}|\nabla\zeta|^2.
\label{estimate-zeta}\eee
We observe that the term $-\lambda\left(\Delta\zeta\nabla\phi,{\bf
w}\right)$ in \eqref{estimate-w} cancels out $\lambda\left(({\bf
w}\cdot\nabla)\phi,\Delta\zeta\right)$ in \eqref{estimate-zeta},
which is originated from the energy law \eqref{basic energy law}.
This is very important, otherwise we cannot control the two terms.

Next we estimate the other terms on the right-hand side of
\eqref{estimate-w}-\eqref{estimate-zeta}, using
\eqref{theta-1}-\eqref{theta-2}. \bee \left|\nu\left(\nabla{\bf
u},\nabla({\bf v}-\bm\theta)\right)\right|&\leq& \nu|\nabla{\bf
u}||\nabla{\bf v}|+\nu|\nabla{\bf
u}|_{L^2(\Gamma_\delta)}|\nabla\bm\theta|_{L^2(\Gamma_\delta)}\non\\
&\leq&\frac{\nu}{4}|\nabla{\bf u}|^2+\nu|\nabla{\bf
v}|^2+C\nu\delta^{-\frac{1}{2}}|\nabla{\bf
u}|_{L^2(\Gamma_\delta)},\\
\left|\left(({\bf w}\cdot\nabla){\bf v}, {\bf w}\right)\right|&\leq&
|\nabla{\bf v}|_{L^\infty}|{\bf w}|^2,\\
\left|\left((\bm\theta\cdot\nabla){\bf v},{\bf
w}\right)\right|&\leq&|\bm\theta||\nabla{\bf v}|_{L^\infty}|{\bf
w}|\leq C\delta^{\frac{1}{2}}|\nabla{\bf v}|_{L^\infty}|{\bf
w}|\non\\
&\leq&C|{\bf w}|^2+\delta|\nabla{\bf v}|_{L^\infty}^2,\eee
\ube \left|\left(({\bf u}\cdot\nabla)\bm\theta,{\bf
u}\right)\right|=\left|\int_{\Gamma_\delta}\left(\frac{1}{\rho}{\bf
u}\cdot\rho^2\nabla\right)\bm\theta\frac{1}{\rho}{\bf
u}dx\right|\leq\left|\frac{1}{\rho}{\bf
u}\right|^2_{L^2(\Gamma_\delta)}\left|\rho^2\nabla\bm\theta\right|_{L^\infty(\Gamma_\delta)}
\leq C\delta|\nabla{\bf u}|^2_{L^2(\Gamma_\delta)}, \uee
\ube \left|\left(({\bf u}\cdot\nabla)\bm\theta,{\bf
v}\right)\right|=\left|\int_{\Gamma_\delta}\left(\frac{1}{\rho}{\bf
u}\cdot\rho\nabla\right)\bm\theta{\bf v}dx\right|
\leq\left|\frac{1}{\rho}{\bf
u}\right|_{L^2(\Gamma_\delta)}\left|\rho\nabla\bm\theta\right||{\bf
v}|_{L^\infty(\Gamma_\delta)}\leq C\delta^{\frac{1}{2}}|\nabla{\bf
u}|_{L^2(\Gamma_\delta)}|{\bf v}|_{L^\infty(\Gamma_\delta)}, \uee
where we used \eqref{theta-2} and the well-known inequality of
Hardy-Littlewood (since ${\bf u}\in H^1_0(\Omega)$).
\be |(\bm\theta_t,{\bf w})|\leq|\bm\theta_t||{\bf w}|\leq
C\delta^{\frac{1}{2}}|{\bf w}|\leq C|{\bf w}|^2+C\delta,\ee \bee
\left|\lambda\left(\Delta\psi\nabla\zeta,{\bf
w}\right)\right|&\leq&\lambda|\Delta\psi|_{L^4}|\nabla\zeta|_{L^4}|{\bf w}|\non\\
&\leq&\left\{ \begin{array}{r@{\quad }l}
C\lambda|\Delta\psi|_{L^4}|\nabla\zeta|^{\frac{1}{2}}
|\Delta\zeta|^{\frac{1}{2}}|{\bf w}| & \quad (\mbox{if}\ n=2)\\
C\lambda|\Delta\psi|_{L^4}|\nabla\zeta|^{\frac{1}{4}}
|\Delta\zeta|^{\frac{3}{4}}|{\bf w}| & \quad (\mbox{if}\ n=3)
\end{array}\right.
\non\\
&\leq&\frac{\lambda\gamma}{8}|\Delta\zeta|^2+C|\nabla\zeta|^2+C|\Delta\psi|_{L^4}^2|{\bf
w}|^2, \eee
\bee \left|\lambda({\bf
w}\cdot\nabla\phi,\zeta)\right|&\leq&\lambda|{\bf
w}||\nabla\phi||\zeta|_{L^\infty}\non\\
&\leq& \left\{ \begin{array}{r@{\quad }l} C|{\bf w}||\nabla\phi||\zeta|^{\frac{1}{2}}|\Delta\zeta|^{\frac{1}{2}} & \quad (\mbox{if}\ n=2)\\
C|{\bf w}||\nabla\phi||\zeta|^{\frac{1}{4}}|\Delta\zeta|^{\frac{3}{4}} & \quad (\mbox{if}\ n=3)\\
\end{array}\right.
\non\\
&\leq& \left\{ \begin{array}{r@{\quad }l}
\frac{\lambda\gamma}{8}|\Delta\zeta|^2+C|{\bf w}|^2+C|\nabla\phi|^4|\zeta|^2 & \quad (\mbox{if}\ n=2)\\
\frac{\lambda\gamma}{8}|\Delta\zeta|^2+C|{\bf
w}|^2+C|\nabla\phi|^8|\zeta|^2&\quad(\mbox{if}\ n=3)
\end{array}\right.
\eee
\bee\left|\lambda({\bm\theta\cdot\nabla}\phi,\zeta)\right|&\leq&\lambda|\bm\theta||\nabla\phi||\zeta|_{L^\infty}\non\\
&\leq&\left\{ \begin{array}{r@{\quad }l}
C\delta^{\frac{1}{2}}|\nabla\phi||\zeta|^{\frac{1}{2}}|\Delta\zeta|^{\frac{1}{2}}
& \quad (\mbox{if}\ n=2)  \\
C\delta^{\frac{1}{2}}|\nabla\phi||\zeta|^{\frac{1}{4}}|\Delta\zeta|^{\frac{3}{4}}
& \quad (\mbox{if}\ n=3)
\end{array}\right.\non\\
&\leq&\frac{\lambda\gamma}{8}|\Delta\zeta|^2+C\delta|\nabla\phi|^2+C|\zeta|^2,
\eee
\bee \left|\lambda\left(({\bf
v}\cdot\nabla)\zeta,\Delta\zeta\right)\right|&\leq&\lambda|{\bf
v}|_{L^4}|\nabla\zeta|_{L^4}|\Delta\zeta|\non\\
&\leq&\left\{\begin{array}{r@{\quad }l}C\lambda|{\bf
v}|_{L^4}|\nabla\zeta|^{\frac{1}{2}}
|\Delta\zeta|^{\frac{3}{2}} & \quad (\mbox{if}\ n=2)\\
C\lambda|{\bf v}|_{L^4}|\nabla\zeta|^{\frac{1}{4}}
|\Delta\zeta|^{\frac{7}{4}} & \quad (\mbox{if}\ n=3)
\end{array}\right.
\non\\
&\leq&\left\{\begin{array}{r@{\quad }l}
\frac{\lambda\gamma}{8}|\Delta\zeta|^2+C|{\bf
v}|_{L^4}^4|\nabla\zeta|^2, & \quad (\mbox{if}\ n=2)\\
\frac{\lambda\gamma}{8}|\Delta\zeta|^2+C|{\bf
v}|_{L^4}^8|\nabla\zeta|^2, & \quad (\mbox{if}\ n=3)\\
\end{array}\right.
\eee
\bee
\left|\lambda\left((\bm\theta\cdot\nabla)\phi,\Delta\zeta\right)\right|&\leq&\lambda|\bm\theta|_{L^4(\Gamma_\delta)}
|\nabla\phi|_{L^4(\Gamma_\delta)}|\Delta\zeta|\leq\lambda\delta^{\frac{1}{4}}
|\nabla\phi|_{L^4(\Gamma_\delta)}|\Delta\zeta|\non\\
&\leq&\frac{\lambda\gamma}{8}|\Delta\zeta|^2+C\delta^{\frac{1}{2}}|\nabla\phi|_{L^4(\Gamma_\delta)}^2,\eee
\bee
\frac{4\lambda\gamma}{\varepsilon^2}\left|\left(\zeta(\phi^2+\phi\psi+\psi^2),\Delta\zeta\right)\right|
&\leq&\frac{4\lambda\gamma}{\varepsilon^2}|\Delta\zeta||\zeta|_{L^6}|\phi^2+\phi\psi+\psi^2|_{L^3}\non\\
&\leq&C|\Delta\zeta||\zeta|_{H^1}\left(|\phi|_{L^6}^2+|\psi|_{L^6}^2\right)\non\\
&\leq&C|\Delta\zeta||\zeta|_{H^1}\left(|\phi|_{H^1}^2+|\psi|_{H^1}^2\right)\non\\
&\leq&\frac{\lambda\gamma}{8}|\Delta\zeta|^2+C\left(|\phi|_{H^1}^4+|\psi|_{H^1}^4\right)|\zeta|_{H^1}^2.
\eee

Adding \eqref{estimate-w}, \eqref{estimate-zeta1} and
\eqref{estimate-zeta}, then using the estimates above and taking
$\delta=c\nu$, we finally obtain \bee
\frac{1}{2}\frac{d}{dt}\left(|{\bf
w}|^2(t)+\lambda|\zeta|^2_1(t)\right)&\leq& \left(C+|\nabla{\bf
v}|_{L^\infty}+C|\Delta\psi|^2_{L^4}\right)|{\bf
w}|^2+C\left(1+|{\bf
v}|^8_{L^4}+|\phi|^8_1+|\psi|^4_1\right)|\zeta|^2_1\non\\
&+&\nu\left(C+|\nabla{\bf v}|^2+|\nabla{\bf
v}|^2_{L^\infty}+C|\nabla\phi|^2\right)+C\nu^\frac{1}{2}\left(1+|{\bf
v}|_{L^\infty(\Gamma_\delta)}\right)|\nabla{\bf u}|_{L^2(\Gamma_\delta)}   \non\\
&+&C\nu|\nabla{\bf
u}|^2_{L^2(\Gamma_\delta)}+C\nu^\frac{1}{2}|\nabla\phi|^2_{L^4(\Gamma_\delta)}.\label{all-estimates}
\eee
Denote \beee Y(t)&=&|{\bf w}|^2+\lambda|\zeta|^2_1,\\
a(t)&=&\nu\left(C+|\nabla{\bf
v}|^2+|\nabla{\bf v}|^2_{L^\infty}+C|\nabla\phi|^2\right),\\
b(t)&=&C\nu^\frac{1}{2}\left(1+|{\bf v}|_{L^\infty(\Gamma_\delta)}
\right)|\nabla{\bf u}|_{L^2(\Gamma_\delta)}+C\nu|\nabla{\bf
u}|^2_{L^2(\Gamma_\delta)}+C\nu^\frac{1}{2}|\nabla\phi|^2_{L^4(\Gamma_\delta)}.\eeee
Since $({\bf v},\psi)$ is a smooth local solution on $[0,T_*]$ and
$\phi$ is bounded in $L^\infty(0,T_*;H^1(\Omega))$, the inequality
\eqref{all-estimates} becomes \be\frac{d Y(t)}{dt}\leq
CY(t)+a(t)+b(t),\ee where $a(t)\leq C\nu$ and
$Y(0)=|\bm\theta|^2(0)$. Recall \eqref{theta-1} that
$|\bm\theta|(t)\leq C\nu^{1/2}$. Hence, \bee \sup_{0\leq t\leq
T_*}Y(t)&\leq& e^{CT_*}\left(|\bm\theta|^2(0)+C\nu
T_*+\int^{T_*}_0b(t)dt\right)\non\\
&\leq&e^{CT_*}\left(C\nu+C\nu T_*+\int^{T_*}_0b(t)dt\right),
\label{Y}\eee where \bee
\int^{T_*}_0b(t)dt&\leq&C{T_*}^{\frac{1}{2}}\left(\nu\int^{T_*}_0|\nabla{\bf
u}|^2_{L^2(\Gamma_\delta)}dt\right)^{\frac{1}{2}}+C\nu\int^{T_*}_0|\nabla{\bf
u}|^2_{L^2(\Gamma_\delta)}dt
+C\nu^{\frac{1}{2}}\int^{T_*}_0|\nabla\phi|^2_{L^4(\Gamma_\delta)}dt\non\\
&\leq&C{T_*}^{\frac{1}{2}}\left(\nu\int^{T_*}_0|\nabla{\bf
u}|^2_{L^2(\Gamma_\delta)}dt\right)^{\frac{1}{2}}+C\nu\int^{T_*}_0|\nabla{\bf
u}|^2_{L^2(\Gamma_\delta)}dt+C\nu^{\frac{1}{2}}, \label{bt}\eee
where we used \eqref{bound} and the fact $H^1(\Omega)\subset
L^4(\Omega)$.
Assume \be \nu\int^{T_*}_0|\nabla{\bf
u}|^2_{L^2(\Gamma_\delta)}dt\rightarrow 0 \ \mbox{as}\, \,
\nu\rightarrow 0.\label{11}\ee Then combining \eqref{Y} and
\eqref{bt}, we have \be \sup_{0\leq t\leq {T_*}}Y(t)\rightarrow 0 \
\mbox{as}\ \nu\rightarrow 0,\label{Yt}\ee which equals to
 \be \sup_{0\leq t\leq {T_*}}\left(|{\bf u}-{\bf v}|^2(t)+\lambda|\phi-\psi|_{H^1}^2(t)\right)\leq \sup_{0\leq t\leq
T_*}Y(t)+|\bm\theta|^2(t)\rightarrow 0 \ \mbox{as}\ \nu\rightarrow
0.\ee

In 2D case, if $({\bf u},\phi)$ is a global strong solution to
\eqref{Navier-Stokes}-\eqref{BC1}, then $|\nabla{\bf u}|^2(t)$ is
uniformly bounded in $[0,\infty)$, thus the assumption \eqref{11} is
automatically satisfied. \eqref{bt} is simplified to \be
\int^{T_*}_0b(t)dt\leq C\nu^{1/2}\,\,\mbox{for}\,\,
\nu\in(0,1).\label{12}\ee Then, using \eqref{12}, \eqref{Y} becomes
\ube \sup_{0\leq t\leq T_*}Y(t)\leq e^{CT_*}\left(C\nu+C\nu
T_*+C\nu^{1/2}\right)\leq C\nu^{1/2} \,\,\mbox{for}\,\, \nu\in(0,1).
\uee Therefore, \be \sup_{0\leq t\leq {T_*}}\left(|{\bf u}-{\bf
v}|^2(t)+\lambda|\phi-\psi|_{H^1}^2(t)\right) \leq \sup_{0\leq t\leq
T_*}Y(t)+|\bm\theta|^2(t)\leq C\nu^{1/2} \,\,\mbox{for}\,\,
\nu\in(0,1). \ee
The proof is complete.

\section{Related Models and Results}
Our system \eqref{Navier-Stokes}-\eqref{Allen Cahn} is closely
related to liquid crystal model, Magnetohydrodynamics (MHD)
equations, and viscoelastic system (see \cite{XZL10} and the
reference therein). In this section, we apply our previous approach
to these models.

\subsection{Liquid Crystal}
If taking $\phi$ as a vector, say ${\bf d}$, system
\eqref{Navier-Stokes}-\eqref{Allen Cahn} can model the motion of
liquid crystal flows: \bee {\bf u}_t+({\bf u}\cdot\nabla){\bf u} &=&
-\nabla p+\nu\Delta{\bf
u}-\lambda\nabla\cdot(\nabla{\bf d}\otimes\nabla{\bf d}), \label{liquid crystal equ 1}\\
\nabla\cdot {\bf u} &=& 0,\\
{\bf d}_t+({\bf u}\cdot\nabla){\bf d} &=& \gamma(\Delta{\bf
d}-f({\bf d})), \label{liquid crystal equ 3} \eee where ${\bf d}$
represents the director field of liquid crystal molecules.
%and boundary conditions:
% \be
% {\bf u}(x,t)=0,\ \ {\bf d}(x,t)= {\bf d}_0(x), \ \ (x,t)\in
%\partial\Omega\times(0,+\infty).
% \ee

We consider the above system in 2D in an infinite channel
$\mathbb{R}\times(0,1)$. We assume the velocity field ${\bf u}$ and
director field ${\bf d}$ of liquid crystals are periodic with period
$2\pi$ in the horizontal ($x$) direction. Set
$\Omega=(0,2\pi)\times(0,1)$. ${\bf u}$ vanishes and the liquid
crystal molecules align in a fixed direction ${\bf d}^*$ at the
boundary of the channel (i.e. at $y=0,1$), \be{\bf u}=0,\ \ {\bf
d}={\bf d}^*, \ \ \mbox{on}\ (0,2\pi)\times\{0,1\}. \ee
\eqref{liquid crystal equ 1}-\eqref{liquid crystal equ 3} is subject
to initial data
 \be
 {\bf u}|_{t=0}={\bf u}_0(x,y),\quad {\bf d}|_{t=0}={\bf d}_0(x,y),\
\mbox{in}\ \Omega \label{BC3}.
 \ee

The above system is a simplified version of the Ericksen-Leslie
model for the hydrodynamics of liquid crystals (cf.
\cite{Ericksen61,Ericksen87,HardtKinder87,Leslie68}). Following
exactly the same arguments as in \cite{LL95}, one obtains the
well-posedness results of system \eqref{liquid crystal equ
1}-\eqref{BC3}.

The corresponding invisid liquid crystal model is \bee {\bf
u}^{(0)}_t+({\bf
u}^{(0)}\cdot\nabla){\bf u}^{(0)} &=& -\nabla P-\lambda\nabla\cdot(\nabla{\bf d}^{(0)}\otimes\nabla{\bf d}^{(0)}), \label{invisid u}\\
\nabla\cdot {\bf u}^{(0)} &=& 0,\\
{\bf d}^{(0)}_t+({\bf u}^{(0)}\cdot\nabla){\bf d}^{(0)} &=&
\gamma(\Delta{\bf d}^{(0)}-f({\bf d}^{(0)})), \label{invisid d}\\
{\bf u}^{(0)},{\bf d}^{(0)}\, \mbox{are}\, 2\pi\mbox{-periodic in}\,
x, && {\bf u}^{(0)}\cdot{\bf n}=0, \quad {\bf d}^{(0)}={\bf d}^*, \
\
\mbox{on}\ (0,2\pi)\times\{0,1\},\\
{\bf u}^{(0)}|_{t=0}={\bf u}_0(x,y), && {\bf d}^{(0)}|_{t=0}={\bf
d}_0(x,y),\ \mbox{in}\ \Omega. \label{BC-d}\eee
Following almost the same proof as that in Theorem \ref{thm3} and
Theorem \ref{thm2}, we have the follows. \bt Let $s\geq 3$ be an
integer. Assume ${\bf u}_0\in X_s$, ${\bf d}_0\in H^{s+1}$. Then
there exists a constant $T_*>0$ depending on $|{\bf u}_0|_s$ and
$|{\bf d}_0|_{s+1}$, such that system \eqref{invisid u}-\eqref{BC-d}
has a smooth solution $({\bf u}^{(0)}, P, {\bf d}^{(0)})$ in
$[0,T_*]$ satisfying \ube {\bf u}^{(0)}\in L^\infty(0,T_*;X_s),\quad
P\in L^\infty(0,T_*; H^{s+1}(\Omega)),\quad {\bf d}^{(0)}\in
L^\infty(0,T_*;H^{s+1}). \uee Moveover, $\left({\bf u}^{(0)}, {\bf
d}^{(0)}\right)$ is unique. \et
\bt Let $({\bf u},{\bf d})$ be a global strong solution to
\eqref{liquid crystal equ 1}-\eqref{BC3}, and $({\bf u}^{(0)},{\bf
d}^{(0)})$ a local smooth solution to \eqref{invisid u}-\eqref{BC-d}
on $\Omega\times[0,T_*]$. Then \be|{\bf u}(t)-{\bf
u}^{(0)}(t)|^2+\lambda|{\bf d}(t)-{\bf d}^{(0)}(t)|_1^2\rightarrow 0
\ \mbox{as}\ \nu\rightarrow 0, \quad \forall
t\in[0,T_*].\label{L-C-conv}\ee The convergence rate in
\eqref{L-C-conv} is $c\nu^{1/2}$ for small $\nu\in(0,1)$. \et

\subsection{The MHD equations}
The nondimensional form of viscous incompressible MHD equations is
(cf. \cite{Temam-MHD1983}) \bee {\bf u}_t+({\bf u}\cdot\nabla){\bf
u}-S({\bf B}\cdot\nabla){\bf B}+\nabla
p&=&\frac{1}{Re}\Delta{\bf u},\\
{\bf B}_t+({\bf u}\cdot\nabla){\bf B}-({\bf
B}\cdot\nabla){\bf u} &=& \frac{1}{Rm}\Delta{\bf B},\\
\nabla\cdot{\bf u}=0, \  \ \nabla\cdot{\bf B}=0,&& \label{Incom of
magnetic field} \eee where ${\bf u}$, $p$ and ${\bf B}$ denote the
fluid velocity, its pressure and the magnetic field, respectively.
Let $L_*, U_*, B_*, \rho_*$ be the characteristic values for
lengths, velocities, magnetic fields and densities of the fluid. The
Reynolds number $R_e=\frac{L_*U_*}{\nu}$ where $\nu$ is the
kinematic viscosity; the magnetic Reynolds number
$R_m=L_*U_*\sigma\mu$ where $\mu$ is the magnetic permeability and
$\sigma$ the resistivity of the fluid;
$S=\frac{M^2}{R_eR_m}\left(=\frac{B_*^2}{\mu\rho_*U_*^2}\right)$,
where $M$ is the Harmann number.

Due to $\nabla\cdot{\bf B}=0$, we know in $2D$ case there exists a
scalar function $\phi$ such that
\[\bf{B}=\nabla^{T}\phi=\left(
                              \begin{array}{c}
                                -\partial_{2} \phi \\
                                 \partial_{1} \phi \\
                              \end{array}
                            \right),
  \]Hence the above  system in $2D$ is
equivalent to the following equations \bee {\bf u}_t+({\bf
u}\cdot\nabla){\bf u}+\nabla p &=& \frac{1}{R_e}\Delta{\bf
u}-S\nabla\cdot(\nabla\phi\otimes\nabla\phi),
\label{equation 1 for MHD}\\
\nabla\cdot {\bf u} &=& 0,  \\
\phi_t+({\bf u}\cdot\nabla)\phi &=& \frac{1}{R_m}\Delta\phi.
\label{equation 3 for MHD} \eee Compared with the system
\eqref{Navier-Stokes}-\eqref{Allen Cahn}, the nonlinear term
$f(\phi)$ vanishes in \eqref{equation 3 for MHD}.

Let $\nu\rightarrow 0$ and other parameters remain constants, we
have $R_e\rightarrow +\infty$ and $S, R_m$ are constants. System
\eqref{equation 1 for MHD}-\eqref{equation 3 for MHD} formally
becomes \bee {\bf v}_t+({\bf v}\cdot\nabla){\bf v}+\nabla P
&=& -S\nabla\cdot(\nabla\psi\otimes\nabla\psi),\label{7}\\
\nabla\cdot {\bf v} &=& 0,  \\
\psi_t+({\bf v}\cdot\nabla)\psi &=&
\frac{1}{R_m}\Delta\psi.\label{8} \eee

We consider the vanishing viscosity limit of system \eqref{equation
1 for MHD}-\eqref{equation 3 for MHD} in the entire plane
$\mathbb{R}^2$ as $R_e\rightarrow +\infty$. Omitting the terms with
boundary layer function ${\bm\theta}$ in the proof of Theorem
\ref{thm2}, and using the same estimates, we easily get similar
result as follows.
\bt Assume $\Omega=\mathbb{R}^2$. Let $({\bf u},\phi)$ be a global
strong solution to \eqref{equation 1 for MHD}-\eqref{equation 3 for
MHD}, and $({\bf v},\psi)$ a local smooth solution to
\eqref{7}-\eqref{8} on $\mathbb{R}^2\times[0,T_*]$. $({\bf u},\phi)$
and $({\bf v},\psi)$ are equipped with the same initial data. Then
\be|{\bf u}(t)-{\bf v}(t)|^2+\lambda|\phi(t)-\psi(t)|_1^2\rightarrow
0 \ \mbox{as}\ R_e\rightarrow +\infty, \quad \forall
t\in[0,T_*].\label{MHD-con}\ee The convergence rate in
\eqref{MHD-con} is $cR_e^{-1}$.

\et
\br In the literature, inviscid limit for MHD equations has been
studied as both $R_e$ and $R_m$ go to infinity, see
\cite{Wu-MHD,MHD-slip2009}. But in some industrial cases (see for
example \cite{FEM-MHD2000}), $R_e\approx 10^5, R_m\approx 10^{-1},
S\approx 1$. Hence it is also meaningful to consider the inviscid
limit for MHD equations as only $R_e\rightarrow +\infty$. \er
%  &\leq& \left( \right)  \left| \right|  |\psi|_{L^\infty}
% \left(\right)  {\bf w}  \left| \right|   L^\infty(0,T_*;X_s)

\subsection{Viscoelasticity}
We consider the following system describing the motion of
incompressible viscoelastic fluids (cf.
\cite{LLZhang05,ZhengLZhou08}): \bee \label{O-Bmodel}
{\bf u}_t+{\bf u}\cdot\nabla{\bf u}+\nabla p&=&\nu\Delta {\bf u}+\lambda\nabla\cdot(\frac{\partial W(\mathcal{F})}{\partial\mathcal{F}}\mathcal{F}^{T}), \\
\nabla\cdot {\bf u}&=&0, \\
\mathcal{F}_t+{\bf u}\cdot\nabla F&=&\nabla {\bf
u}\mathcal{F}.\label{F} \eee Here, ${\bf u}$ represents the velocity
field of materials, $p$ the pressure, $\nu$ the kinetic viscosity
constant, and $\lambda$ a parameter representing the competition
between kinetic energy and elastic energy. $\mathcal{F}$ is the
deformation tensor, and $W(\mathcal{F})$ the elastic energy
functional.

\br The dynamic of viscoelastic fluids can be described by the flow
map. Let $X$ be the original labeling (Lagrangian) coordinate of the
particle, and $x$ the current (Eulerian) coordinate. For a given
velocity ${\bf u}(x,t)$, the flow map (particle trajectory) $x(X,t)$
is defined by the ODE :
\[x_t={\bf u}(x(X, t), t), \;\; x(X, 0)=X. \] The deformation tensor
$\tilde{\mathcal{F}}(X,t)$ is defined as
$\tilde{\mathcal{F}}(X,t)=\frac{\partial x}{\partial X}$, where we
use the notation $\tilde{\mathcal{F}}_{ij}=\frac{\partial
x_i}{\partial X_j}$. In the Eulerian coordinate, we can define
$\mathcal{F}(x, t)=\mathcal{F}(x(X,t),t)=\tilde{\mathcal{F}}(X, t)$.
The equation \eqref{F} is deduced by the chain rule. The equation
\eqref{O-Bmodel} can be derived by the energetic variational
approach, see \cite{LLZhang05,XZL10} for intrinsic mechanics in
detail.\er

Taking divergence of both sides of \eqref{F} and using
$\nabla\cdot{\bf u}=0$, we get the transport equation for
$\nabla\cdot\mathcal{F}$ as \be (\nabla\cdot\mathcal{F})_t+{\bf
u}\cdot\nabla(\nabla\cdot\mathcal{F})=0.\label{divergence of
deformation tensor F} \ee Since $\mathcal{F}$ is the identity matrix
at $t=0$, $\nabla\cdot\mathcal{F}|_{t=0}=0$, then \eqref{divergence
of deformation tensor F} tells that \be \label{zero divergence of
deformation tensor F} \nabla\cdot\mathcal{F}=0,\ \ \forall t \geq 0
. \ee In 2D, this implies there exists a vector function
$\bm\phi=(\phi_1,\phi_2)$, such that
\be\mathcal{F}=\nabla^{T}\bm\phi=\left(
                              \begin{array}{c}
                                -\partial_{2} \bm\phi \\
                                 \partial_{1} \bm\phi \\
                              \end{array}
                            \right).  \label{F-2D}\ee
 If we take
$W(\mathcal{F})$ as the Hookean elasticity energy, i.e.,
$W(\mathcal{F})=\frac{1}{2}tr(\mathcal{F}\mathcal{F}^T)$, and use
\eqref{F-2D}, then the system \eqref{O-Bmodel}-\eqref{F} is
equivalent to \bee{\bf u}_t+({\bf u}\cdot\nabla){\bf u} &=& -\nabla
q+\nu\Delta{\bf
u}-\lambda\nabla\cdot(\nabla\bm\phi\otimes\nabla\bm\phi), \label{Navier Stokes in Oldroyd-B}\\
\nabla\cdot {\bf u} &=& 0,  \label{imcompressibility in Oldroyd-B}\\
\bm\phi_t+({\bf u}\cdot\nabla)\bm\phi &=& 0, \label{trans equ in
Oldroyd} \eee where the basic energy law is \bee
\frac{1}{2}\frac{d}{dt}\left(|{\bf
u}|^2+\lambda|\nabla\bm\phi|^2\right)=-\nu|\nabla {\bf u}|^2. \eee

Due to the system has partial dissipation, it is difficult to obtain
global classical solution in general. Consider the system
\eqref{Navier Stokes in Oldroyd-B}-\eqref{trans equ in Oldroyd} in
the entire plane $\mathbb{R}^2$. The existence of global smooth
solutions near equilibrium has been proved in \cite{LLZhang05}.
Moreover, the authors proved the local existence of a classical
solution near equilibrium to the inviscid system in $\mathbb{R}^2$:
\beee{\bf u}_t+({\bf u}\cdot\nabla){\bf u} &=& -\nabla q-\lambda\nabla\cdot(\nabla\bm\phi\otimes\nabla\bm\phi),\\
\nabla\cdot {\bf u} &=& 0,  \\
\bm\phi_t+({\bf u}\cdot\nabla)\bm\phi &=& 0, \eeee by rewritting the
above system as a quasilinear first order system.

Due to the lack of dissipation term in the equation for $\phi$, we
cannot show the validity of vanishing viscosity limit for system
\eqref{Navier Stokes in Oldroyd-B}-\eqref{trans equ in Oldroyd}, by
similar proof to that of Theorem \ref{thm2}. But we believe that the
inviscid limit holds for solutions near equilibrium. We are going to
investigate this problem in forthcoming work.

\noindent\textbf{Acknowledgements.} The research of L. Zhao and H.
Huang is partially supported by National Natural Science Foundation
of China under Grant No. 10971142.

\end{document}